\documentclass[12pt]{article}
\usepackage{amsmath, amssymb, amsthm, amscd}

\numberwithin{equation}{section}

\theoremstyle{plain}
\newtheorem{theorem}{Theorem}
\newtheorem{proposition}{Proposition}
\newtheorem{lemma}{Lemma}
\newtheorem{corollary}{Corollary}
\newtheorem{remark}{Remark}
\theoremstyle{definition}
\newtheorem{definition}{Definition}

\def\sett{\sqrt{ \Omega_s}}

\def\C{{\mathbb C}}
\def\R{{\mathbb R}}

\def\K{{\mathfrak K}}
\def\Hh{{\cal H}^{\rm{\scriptscriptstyle H}}}
\def\dh{\delta^{\rm{\scriptscriptstyle H}}}
\def\db{\delta^{\rm{\scriptscriptstyle BKS}}}
\def\bks{\rm{\scriptscriptstyle BKS}}
\def\Hpr{{\cal H}^{\rm{\scriptscriptstyle prQ}}}

\def\dpr{\delta^{\rm{\scriptscriptstyle prQ}}}
\def\Hq{{\cal H}^{\rm{\scriptscriptstyle Q}}}

\def\dq{\delta^{\rm{\scriptscriptstyle Q}}}

\def\Qu{\rm{\scriptscriptstyle Q}}
\def\Pr{\rm{\scriptscriptstyle prQ}}

\def\bi{\begin{itemize}}
\def\ei{\end{itemize}}
\def\la{\label}

\newcommand{\be}{\begin{equation}}
\newcommand{\ee}{\end{equation}}
\newcommand{\ba}{\begin{eqnarray}}
\newcommand{\ea}{\end{eqnarray}}

\title{On the BKS Pairing for K\"ahler Quantizations of the Cotangent Bundle 
of a Lie Group}
\author{Carlos Florentino$^\dagger$, Pedro Matias$^\ddagger$,
Jos\'e Mour\~ao$^{\dagger}$ \\ and Jo\~ao P.
Nunes$^\dagger$}

\begin{document}

\maketitle

\begin{abstract}

A natural one-parameter family of 
K\"ahler quantizations of the cotangent bundle
$T^*K$ of a compact Lie group $K$, taking into account the half-form correction,
 was studied in \cite{FMMN}.
In the present paper, it is shown that the associated Blattner-Kostant-Sternberg (BKS) 
pairing map is unitary and coincides with the parallel transport of 
the quantum connection introduced in 
our previous work, from the point of view of \cite{AdPW}.
The BKS pairing map is a composition of (unitary) 
coherent state transforms of $K$, introduced in \cite{Ha1}. 
Continuity of the Hermitian structure on the quantum bundle,
in the limit when one of the K\"ahler polarizations 
degenerates to the vertical real polarization,
leads to the unitarity of the corresponding BKS 
pairing map. 
This is in agreement with 
the 
unitarity up to scaling (with respect to a rescaled inner product) 
of this pairing map, 
established by Hall.

\end{abstract}

\newpage

\tableofcontents

\section{Introduction}\label{s1}

Let $K$ be a compact, connected Lie group of dimension $n$ and let $T^*K$ be its cotangent bundle.
We start by recalling  some aspects of \cite{FMMN}
where, in connection with work of Hall in \cite{Ha3}, the
geometric quantization of $T^*K$ was studied using 
a natural one-parameter family of K\"ahler structures.
These K\"ahler structures are induced on $T^*K$ via the following
natural identifications of $T^*K$ with the complexified group $K_\C$.
Consider, for any real parameter $s>0$, the 
diffeomorphisms
\ba
\label{difeo}
\psi_s: T^*K & \rightarrow & K_\C \nonumber \\
(x,Y)&\mapsto & \psi_s(x,Y)=xe^{isY}\ .
\ea
Here, $x\in K$, $Y\in \K\equiv{\rm Lie}(K)$, and we identify $T^*K$ with $K\times\K^*$ 
using left invariant forms
and then with $K\times\K$ by means of a fixed $Ad$-invariant inner product
on $\K$. 
The diffeomorphisms $\psi_s$ endow $T^*K$ with a family of complex structures $J_s$ and 
one can check that, together with the canonical symplectic structure $\omega$ on $T^*K$,
the pair $(\omega,J_s)$ defines a K\"ahler structure on $T^*K$ for every $s \in \R_+$ \cite{FMMN}.
This family includes the K\"ahler structure on $T^*K$
considered by Hall in \cite{Ha3}.

In this paper, we consider the Blattner-Kostant-Sternberg (BKS)
pairing between two different K\"ahler quantizations of $T^*K$. 
To describe the results, let
us consider the framework used in \cite{FMMN}.
Let $L$ denote the trivial complex line bundle
on $T^*K$, with trivial Hermitian
structure (its sections are therefore identified with $C^\infty$ functions
on $T^*K$). Following the geometric
quantization program with half-form correction, 
let us introduce the half-form bundle $\delta_s$, which is a square
root of the (trivial) $J_s$-canonical bundle $\kappa_s$ over $T^*K$.
Choosing canonical trivializing $J_s$-holomorphic sections $\Omega_s$ of $\kappa_s=\delta_s^2$ 
and $\sqrt{\Omega_s}$ of $\delta_s$ (we refer to the next section for precise formulas) one
introduces a natural Hermitian structure on $L\otimes\delta_s$ so that,
for a smooth section $\sigma_s$ of the form 
\be
\label{parsec}
\sigma_s= f \sett ,
\ee
with $f\in C^\infty(T^*K)$, one has
\be
\label{hermit}
|\sigma_s|^2:=|f|^2 |\Omega_s| ,
\ee
where
$|\Omega_s|$ is defined by $\overline \Omega_s \wedge \Omega_s = |\Omega_s|^2 b\epsilon $,
$b=(2i)^{n}(-1)^{n(n-1)/2}$ and 
$\epsilon=\frac1{n!}\omega^n$ is the Liouville measure on $T^*K$.

The prequantum Hilbert space $\Hpr_s$, depending on $s$, is then
the (norm completion of the) space of $C^\infty$ sections $\sigma_s$ of
$L\otimes \delta_s$ which are square-integrable with respect
to the Hermitian structure (\ref{hermit}), that is
\begin{equation}
\label{e37}
\langle\sigma_s,\sigma_s\rangle^{\Pr} := \int_{T^*K} |\sigma_s|^2 \epsilon
\quad < +\infty .
\end{equation}

Proceeding with the quantization program, 
the quantum Hilbert space $\Hq_s$ is defined to be the subspace of $\Hpr_s$ consisting
of polarized ($J_s$-holomorphic) sections of $L\otimes \delta_s$. This is
naturally a sub-Hilbert space of $\Hpr_s$.

Both families
of Hilbert spaces can be collected to form the Hilbert prequantum bundle
$\Hpr\to\R_+$ and the quantum bundle $\Hq\to\R_+$, which is naturally
a sub-bundle of $\Hpr$. 

In the spirit of  \cite{AdPW}, we consider a natural connection $\dpr$ on $\Hpr$
and use it to induce, by orthogonal projection, an Hermitian connection 
$\dq$ on $\Hq$. Note that, from
(\ref{parsec}-\ref{e37}), 
$\Hpr$ has a natural global (inner product preserving)  trivializing morphism 
 defined by
\ba
\nonumber
L^2(T^*K,\epsilon) \times \R_+& \to & \Hpr \\
\nonumber
(f,s) & \mapsto &  \frac{f}{\sqrt{|\Omega_s|}} \sett.
\ea
The prequantum connection $\dpr$ on 
$\Hpr$ is defined to be the connection induced via this map from
the trivial bundle $L^2(T^*K,\epsilon) \times \R_+ \to \R_+$ equipped with the 
trivial connection.
Then, the quantum connection $\dq$
is defined to be simply the orthogonal projection of $\dpr$
to the quantum bundle.
This construction is a natural generalization, to the case with the half-form correction and 
with a natural trivializing section of the square root of the canonical bundle, 
of the framework considered in \cite{AdPW}. One of the
advantages of this approach consists in the fact that it ensures  
automatically that the quantum connection is Hermitian.

In the present paper, the framework of \cite{AdPW} is extended in another natural way.
Namely, we consider the Blattner-Kostant-Sternberg (BKS)
pairing between two different fibers of $\Hq$, 
corresponding to two different quantizations of $T^*K$ (see, for instance \cite{Wo}). 
This is the restriction of a Hermitian pairing 
between the corresponding prequantum Hilbert spaces, 
which we call the prequantum BKS pairing, and which
is defined by 
\be
\nonumber
%\label{bks}
\langle \sigma_s,\sigma'_{s'}\rangle^{\bks} = \int_{T^*K} \bar f f' 
\sqrt{\frac{\bar \Omega_s \wedge \Omega_{s'}}{b\epsilon}} \epsilon,
\ee
for $\sigma_s=f\sqrt{\Omega_s}\in\Hpr_s$ and 
$\sigma'_{s'}=f'\sqrt{\Omega_{s'}}\in\Hpr_{s'}$.
This prequantum BKS pairing defines a map
\be
\label{cff}
B^{\Pr}_{ss'}: \Hpr_{s'}\rightarrow \Hpr_s,
\ee
whose infinitesimal version induces a, also natural, connection $\delta$ on $\Hpr$. 
We prove in theorem \ref{bks8}  that $\delta$
coincides with the connection $\dpr$. 
Therefore, the quantum connection $\dq$ will also be induced from the BKS pairing on $\Hq$.
This result shows
that the approach of the infinitesimal BKS pairing coincides with that of
\cite{AdPW} considered in \cite{FMMN}. We note, however, that the
prequantum BKS pairing map does not coincide with the parallel transport 
associated to $\dpr$ (see theorem \ref{bksze}.)

In \cite{Ha3}, Hall considered two polarizations for the quantization
of $T^*K$. Fixing the same inner product as Hall on $\K$, one of these
polarizations, the K\"ahler one, corresponds to the case of $s=1$, and the other,
the vertical real polarization, corresponds, in our framework, 
to letting $s$ go to 0. 
He proved that the BKS pairing map between the Hilbert spaces associated to 
these two quantizations coincides (up to scale) with the Segal-Bargmann
coherent state transform (CST) for $K$ introduced in \cite{Ha1},
provided that one takes into account the
half-form correction.
In theorem \ref{ze99}, in the third section of this paper,
we elaborate on this relation, by showing that the (quantum) 
BKS pairing map
\be
\nonumber
B^{\Qu}_{ss'}: \Hq_{s'}\rightarrow \Hq_s
\ee 
obtained by restriction of (\ref{cff}) to the quantum bundle is, in fact, 
unitary. Moreover, and unlike the prequantum case, 
we show in theorem \ref{bks145} that $B^{\Qu}_{ss'}$ coincides with the parallel
transport of $\dq$ and corresponds also to the CST.
In particular this implies that the BKS pairing maps satisfy
$$
B^{\Qu}_{ss'} \circ B^{\Qu}_{s's''} = B^{\Qu}_{ss''} .
$$

The quantum bundle $\Hq$ can be naturally extended to a bundle 
over $[0,\infty)$ by including a fiber over $s=0$ corresponding
to the space of vertically polarized sections.
Continuity at $s=0$ of the Hermitian structure on the extended 
quantum bundle motivates a particular choice of scale of the inner 
product on the space of vertically polarized sections.
We show that with this choice the BKS pairing map $B^{\Qu}_{s0}$
becomes unitary.
We note that this is in agreement with \cite{Ha3}, where this 
map was shown to be unitary up to scale but for a different 
choice of scale in the inner product.
We study this
fact in section \ref{ss31aq}.

Therefore, we have obtained a class of
examples of quantization where the BKS pairing map, between 
half-form corrected Hilbert spaces corresponding to two different 
K\"ahler quantizations, is 
an isometric isomorphism. We also confirm in a different way the result of
\cite{FMMN} that the parallel transport of the quantum connection coincides
with the one given by the CST. 
We stress that we prove unitarity of the quantum BKS pairing map and not just unitarity
up to scale. Therefore, for the family of polarizations that we consider, 
cotangent bundles of compact Lie groups provide new 
examples in which the half-form correction leads to the unitarity 
of the BKS pairing map (which is the case for Fock spaces, see \cite{Wo}, Ch. 10).

For other related works on the dependence of the quantization on the choice of polarization, the 
unitarity of the BKS pairing map and the relation between vertical and K\"ahler polarizations, 
see also \cite{FY,Hi,KW,Ra,Th}.

\section{The prequantum BKS pairing}\label{s2}
\label{bkspair}

Let $K$ be a compact, connected Lie group of real
dimension $n$, and let us fix once and for all the
$Ad$-invariant inner product $(\cdot\, ,\cdot )$ on its
Lie algebra $\K$, for which the corresponding 
Riemannian measure is the normalized Haar measure on $K$.
Let $\{X_i\}_{i=1}^{n}$ be left-invariant vector fields 
on $K$ forming an orthonormal basis for $\K$.
The cotangent bundle $T^*K \cong K \times \K$ is naturally a symplectic manifold
 with  a canonical
symplectic 2-form  defined by $\omega=-{\rm d}\theta$, where
\be
\nonumber
\theta=\sum_{i=1}^{n}y^{i}w^i.
\ee 
Here,
$(y^1,\ldots ,y^n)$ are the global coordinates on $\K$
corresponding to the basis $\{X_i\}_{i=1}^n$, 
and $\{w^i \}_{i=1}^n$ is the basis of left-invariant
1-forms on $K$ dual to $\{X_i\}_{i=1}^n$, pulled-back to $T^*K$ by the 
canonical projection.
We let
$\epsilon$ denote the Liouville volume form on $T^*K$, given by
\begin{equation}
\nonumber
%\label{e30}
\epsilon=\frac{1}{n!}\omega^{n} \ .
%\hbox{(n times)}
\end{equation}
{} For each $s>0$, let $\psi_s: T^*K \rightarrow K_\C$ denote the diffeomorphism
defined in (\ref{difeo}) by $\psi_s(x,Y)=xe^{isY}$, which induces
on $T^*K$ a K\"ahler structure $(\omega,J_s)$.

\subsection{The prequantum bundle $\Hpr$}

We start by recalling from \cite{FMMN}, some of the formulas that will be needed later on.
Let $\tilde{X_j}$, $j=1,...,n$, be the vector fields on $T^*K$
generating the right action of $K$ 
lifted to $T^*K$ and given by
%\begin{equation}
%\label{lift}
%\tilde{X_j}(x,Y)=(X_j, [Y,X_j]),
%\end{equation}
%so that,
\begin{equation*}
\label{lift}
\psi_{s*}\tilde{X_j}= X_{j,\C},
\end{equation*}
where $X_{j,\C}$ denotes the extension of
$X_j$ from a left-invariant vector field
on $K\subset K_\C$ to  the corresponding  left-invariant vector field on $K_\C$.
Let $\{\tilde w^j\}$ be the basis of left invariant  
1-forms defined by $\tilde w^j (\tilde X_k) = \delta^j_k$ and 
$\tilde w^j (J_s \tilde X_k)=0$, for $j,k = 1,...,n$.  
{}For every $s \in \R_+$,
a frame of left invariant $J_s$-holomorphic  1-forms is then 
\begin{equation*}
%\label{formas}
\left\{\tilde \eta^j_{s} =  \tilde w^j - i J_s \tilde w^j \right\}_{j=1}^n \ ,
\end{equation*}
where $(J_s  w)(X)=w(J_s X)$, for a vector field $X$ and a 1-form $w$
on $T^*K$.

Consider now the $J_s$-canonical bundle $\kappa_s\to T^*K$ whose sections
are $J_s$-holomorphic $n$-forms with natural Hermitian
structure defined as follows.
For a $J_s$-holomorphic
$n$-form $\alpha_s$, let $|\alpha_s|$ be the unique non-negative 
$C^\infty$ function on $T^*K$ such that $\overline \alpha_s \wedge
\alpha_s = |\alpha_s|^2 b\epsilon$, where $b=(2i)^{n}(-1)^{n(n-1)/2}$.
Following \cite{Ha3} we write
\be
\nonumber
%\label{bks1}
|\alpha_s|^2 = \frac{\overline \alpha_s \wedge \alpha_s}{b\epsilon}.
\ee

Given the inner product on $\K$, a canonical trivializing 
$J_s$-holomorphic section of $\kappa_s$ is given by
\begin{equation*}
%\label{e33}
\Omega_s := \tilde \eta_s^1 \wedge \cdots \wedge \tilde \eta_s^n
\end{equation*}
and has norm \cite{FMMN}
\be
\label{normaOmega}
|\Omega_s|^2= s^{n} \eta^2(sY) ,
\ee
where $\eta(Y)$
is the $Ad_K$-invariant function, 
defined for $Y$ in a Cartan subalgebra
by the following product over a set $R^+$ of positive roots of $\K$,
\begin{equation}\label{e25}
\eta (Y)=\prod_{\alpha\in R^+} \frac{\sinh \alpha (Y)}{\alpha (Y)}.
\end{equation}

The following proposition generalizes (\ref{normaOmega}).

\begin{proposition}
\label{bks4}
Let $s,s'>0$. We have 
\be
\nonumber
%\label{bks5}
\bar \Omega_s \wedge \Omega_{s'} = \left( \frac{s+s'}{2}\right)^n 
\eta^2\left(\frac{s+s'}2 Y\right)
b\epsilon = \left| \Omega_{\frac{s+s'}2} \right|^2 b\epsilon . 
\ee
\end{proposition}

\begin{proof}
The result follows by direct computation.
{}From \cite{Ha2} and the definition of $\psi_s$, 
we can write $D\psi_s:TT^*K\to TK_\C$ as 
\begin{equation}
\nonumber
%\label{ap1}
D\psi_s (x,Y) = \left[
\begin{array}{cc}
\cos {\rm ad} sY & \frac{1-\cos {\rm ad} sY}{{\rm ad} Y}\\[2mm]
-\sin {\rm ad} sY & \frac{\sin {\rm ad} sY}{{\rm ad} Y}
\end{array}\right],
\end{equation}
using the $(x,Y)$ coordinate basis on $T^*K$ and the basis
$\{X_{j,\C},JX_{j,\C}\}_{j=1,...,n}$ on $K_\C$,
where $J$ is the complex structure on $K_\C$.
{}From the explicit expressions for the forms 
$\tilde \eta_s^j$ we then find, for $s,s'>0$, 
\be
\nonumber
%\label{ap2}
\bar \Omega_s \wedge \Omega_{s'} = 
\bar {\tilde \eta}_s^1 \wedge \cdots \wedge \bar{\tilde \eta}_s^n 
\wedge \tilde \eta_{s'}^1 \wedge \cdots \wedge \tilde \eta_{s'}^n =
\det \left[
\begin{array}{cc}
\bar M_s & \bar N_s \\
M_{s'} & N_{s'}
\end{array}\right] (-1)^{\frac{n(n-1)}{2}}\epsilon,
\ee
where the endomorphisms $M_s$ and $N_s$ are defined by 
$$
M_s = e^{-i {\rm ad} sY}, \quad N_s = \frac{1-e^{-i {\rm ad} sY}}{{\rm ad}Y}. 
$$

{}From the left invariance of the forms, this determinant can be evaluated 
for $Y$ in the Cartan subalgebra which, after taking care of the 
contribution from the null space of ${\rm ad Y}$, yields
\[
\bar \Omega_s \wedge \Omega_{s'} = \left( \frac{s+s'}{2}\right)^n 
\prod_{\alpha\in R} \frac{(e^{s'\langle\alpha,Y\rangle}-e^{-s\langle \alpha,Y\rangle})}
{(s+s')\langle\alpha,Y\rangle} 
\, b \epsilon,
\]
where the product runs over the set $R$ of all roots of $\K$. 
The result then follows from definition (\ref{e25}).
\end{proof}

As in the introduction, let $\delta_s$ be the $J_s$-holomorphic bundle  
of half-forms on $T^*K$, with trivializing section whose square
is $\Omega_s$. Following
\cite{Ha3}, we will denote this section by $\sqrt{\Omega_s}$. 
Recall that the prequantum Hilbert space $\Hpr$ is the completion
of the space of $C^\infty$-sections of $L\otimes \delta_s$ of finite norm
(\ref{e37}), where $L$ denotes the trivial complex line bundle on $T^*K$, and 
its Hermitian structure $\langle \cdot,\cdot\rangle^{\Pr}$ 
can be written as
\be
\nonumber
%\label{bks2}
\langle \sigma_s,\sigma'_s\rangle^{\Pr} = \int_{T^*K} \bar f f' |\Omega_s| \epsilon =
\int_{T^*K} \bar f f' \sqrt{\frac{\bar \Omega_s\wedge \Omega_s} 
{b\epsilon}} \epsilon,
\ee
for two smooth sections of $L\otimes\delta_s$ written as 
$\sigma_s = f\sqrt{\Omega_s}, \sigma'_s = f' \sqrt{\Omega_s}$, 
with $f,f'\in C^\infty (T^*K)$.

The smooth Hilbert bundle structure on $\Hpr$ and the prequantum connection
$\dpr$ are chosen to be the ones 
compatible with the global trivializing map
\ba
\label{bks6}
L^2(T^*K,\epsilon) \times \R_+ & \to & \Hpr \\
\nonumber
(f,s) & \mapsto & \frac{f}{\sqrt{|\Omega_s|}}\sqrt{\Omega_s},
\ea
and the trivial connection on $L^2(T^*K,\epsilon) \times \R_+ $.
Note that the section of $\delta_s$ given by $\frac{\sett}{\sqrt{|\Omega_s|}}$
has unit norm, and its use would simplify some of the formulas. However,
since this section is not holomorphic, we use $\sett$ as
trivializing section of $\delta_s$, which will be more suited for
describing the polarized sections in section 3.

\subsection{The prequantum BKS pairing and its associated connection on $\Hpr$}

The general procedure, as described for instance in  \cite{Wo}, for defining
a BKS pairing, suggests, in our setting, the definition of a Hermitian
pairing which we call the prequantum BKS pairing, as follows.

\begin{definition}
\label{preqbks}
Let $\sigma_s = f\sqrt{\Omega_s}\in \Hpr_s$ and $\sigma'_{s'} = f'\sqrt{\Omega_{s'}}\in \Hpr_{s'}$. 
Their BKS pairing is defined by 
\be
\label{bks3}
\langle \sigma_s,\sigma'_{s'}\rangle^{\bks} = \int_{T^*K} \bar f f' 
\sqrt{\frac{\bar \Omega_s \wedge \Omega_{s'}}{b\epsilon}} \epsilon.
\ee
\end{definition}

Note that the integral above exists for $\sigma_s,\sigma'_{s'}$
satisfying the conditions for sections of the prequantum bundle in (\ref{e37}), which in this case 
read 
$f\sqrt{|\Omega_s|}, f'\sqrt{|\Omega_{s'}|}\in L^2(T^*K,\epsilon)$. 
This can be readily checked by using Proposition \ref{bks4} to
write (\ref{bks3}) as
\be
\label{bks-1}
\langle \sigma_s,\sigma'_{s'}\rangle^{\bks} = \int_{T^*K} \bar f f' 
\left|\Omega_{\frac{s+s'}2}\right| \epsilon ,
\ee
and using the definition of $\eta(Y)$ in (\ref{e25}) to prove that the smooth function 
\be
\label{ap3}
\phi (s,s',Y) := \frac{\left|\Omega_{\frac{s+s'}2}\right|^2}{|\Omega_s||\Omega_{s'}|} =
\left(\frac{s+s'}{2\sqrt{s s'}}\right)^n  \frac{ \eta^2\left(\frac{s+s'}2 Y\right)}
{\eta(sY)\eta(s'Y)}
\ee
is real, positive and bounded, for fixed $s,s'$.

We remark that, in the case when $s'$ and $s$ coincide,
the prequantum BKS pairing is equal to the Hermitian structure 
$\langle \cdot, \cdot\rangle^{\Pr}$ on $\Hpr_s$.
Also, we note that the pairing (\ref{bks3}) is nondegenerate.

As mentioned above, $\dpr$ was defined 
as the Hermitian connection for which sections $\Hpr$ of the 
form $\sigma_s = \frac{f}{\sqrt{|\Omega_s|}}\sqrt{\Omega_s}$, 
with $f\in L^2(T^*K,\epsilon)$, are parallel. 
It is also natural to consider a connection, $\delta$, on $\Hpr$, induced from the infinitesimal 
prequantum BKS pairing by the formula
\be
\label{bks7}
\langle \sigma_s, \delta_{\frac{\partial}{\partial s}}\sigma'_{s}\rangle^{\Pr} = \left.
\frac{\partial}{\partial s'}\right|_{s'=s} \langle \sigma_s, \sigma'_{s'}\rangle^{\bks}. 
\ee  
Let us prove that these two connections are the same.

\begin{theorem}
\label{bks8}
The connections $\dpr$ and $\delta$ on $\Hpr$ are equal.
\end{theorem}

\begin{proof}
Let $\sigma$ be a smooth section of $\Hpr$, such that 
$\sigma_s = \frac{f}{\sqrt{|\Omega_s|}}\sqrt{\Omega_s} 
\in \Hpr_s$, for $f\in L^2 (T^*K,\epsilon)$, $s>0$.
Since sections of this type can be used to form a global moving 
frame for $\Hpr$, to prove the theorem, it suffices to show that 
$\sigma$ is parallel with respect to the 
connection $\delta$.
Let $\tau_s= g \sqrt{\Omega_s} \in \Hpr_s$. We have, from
(\ref{bks-1}) and (\ref{ap3}),
$$
\langle \tau_s, \delta_{\frac{\partial}{\partial s}}\sigma_{s}\rangle^{\Pr} 
= \left. \frac{\partial}{\partial s'}\right|_{s'=s} \langle \tau_s, \sigma_{s'}\rangle^{\bks} = 
\int_{T^*K} \bar g f |\Omega_s|^{\frac12}\left( \frac{\partial 
\sqrt{\phi}}{\partial s'}\right)_{\!\!s'=s} 
\epsilon.
$$
{}From a straightforward computation using (\ref{ap3}) we obtain
$$
\nonumber
\left. {\frac{\partial \phi (s,s',Y)}{\partial s'}}\right|_{{s' =s}} = 0,
$$
which implies that 
$\delta_{\frac{\partial}{\partial s}}\sigma_{s}=0$, 
as required.
\end{proof}

One immediate consequence of this theorem is that
 $\dpr$ induced from the infinitesimal prequantum BKS pairing is unitary. 
We note that the calculation in the proof of theorem \ref{bks8} relies 
just on the expression for 
$\bar \Omega_s \wedge \Omega_{s'}$
in proposition \ref{bks4}. While in this case this is related 
to the fact that the continuous family of complex polarizations
is generated by the flux of a particular vector field \cite{FMMN}, 
it is not clear to us whether this holds more generally.

In the next section, we will show that the restriction of the prequantum BKS pairing to 
$\Hq$ yields the (unitary) parallel transport associated to $\dq$. 
This, however, does not hold for the prequantum bundle. 
In fact, let us define the prequantum BKS pairing map
$B^{\Pr}_{ss'}: \Hpr_{s'}\rightarrow \Hpr_s$ by 
\be
\label{BKSmap}
\langle \sigma_s,\sigma'_{s'}\rangle^{\bks} = \langle \sigma_s, B^{\Pr}_{ss'}\sigma'_{s'}\rangle^{\Pr},
\ee

\begin{theorem}
\label{chapeu}
The prequantum BKS pairing map (\ref{BKSmap}), $B^{\Pr}_{ss'}: \Hpr_{s'}\rightarrow \Hpr_s$
is given by 
\be
\nonumber
B^{\Pr}_{ss'}\left( f'\frac{\sqrt{\Omega_{s'}}}{\sqrt{|\Omega_{s'}|}}\right)=  
\sqrt{\phi}f'\frac{\sqrt{\Omega_s}}{\sqrt{|\Omega_{s}|}},
\ee
for $f'\frac{\sqrt{\Omega_{s'}}}{\sqrt{|\Omega_{s'}}|}\in \Hpr_{s'}$, where $\phi$ is given 
by (\ref{ap3}).
Consequently, $B^{\Pr}_{ss'}$ 
does not coincide with the parallel transport of the connection $\dpr$. 
\label{bksze}
\end{theorem}

\begin{proof}
Let us write 
\[
B_{ss'}^{\Pr}\left(f'\frac{\sqrt{\Omega_{s'}}}{\sqrt{|\Omega_{s'}|}}\right)
= g\frac{\sqrt{\Omega_s}}{\sqrt{|\Omega_{s}|}},
\]
for some vector $g\frac{\sqrt{\Omega_s}}{\sqrt{|\Omega_{s}|}}\in \Hpr_s$. 
Then, using (\ref{bks-1}) we have for any
$f\sqrt{\Omega_s}\in \Hpr_s$, that 
$$
\int_{T^*K} \bar ff' \frac{\left|\Omega_{\frac{s+s'}{2}}\right|}{\sqrt{|\Omega_{s'}|}} \epsilon = 
\langle f\sqrt{\Omega_s}, B_{ss'}^{\Pr} (f'\frac{\sqrt{\Omega_{s'}}}{\sqrt{|\Omega_{s'}|}})\rangle^{\Pr} =
\int_{T^*K} \bar f g \sqrt{|\Omega_s|}\epsilon.
$$
Since $f\in L^2(T^*K,|\Omega_s|^{\frac12}\epsilon)$ 
is arbitrary, we obtain $g = \sqrt{\phi}f'$, as wanted.

On the other hand, the parallel transport $P_{ss'}: \Hpr_{s'}\rightarrow \Hpr_s$ of $\dpr$ is determined by
(\ref{bks6}) giving
\be
\nonumber
P_{ss'} (f'\frac{\sqrt{\Omega_{s'}}}{\sqrt{|\Omega_{s'}|}}) =  
f' \frac{\sqrt{\Omega_s}}{\sqrt{|\Omega_{s}|}}.
\ee
Clearly, this is not equal to $B_{ss'}^{\Pr}$ since $\phi$ is not identically equal to 1.
\end{proof}

It follows from theorem \ref{chapeu} that $B_{ss'}^{\Pr}$ is not unitary. We will see, however, 
that its restriction to the quantum Hilbert bundle is unitary.

\section{Unitarity of the quantum BKS pairing map and the CST}\label{s3}

Consider the quantum Hilbert bundle $\Hq \rightarrow \R_+$, 
whose fiber over each $s>0$ consists of $J_s$-polarized sections of $L\otimes\delta_s$. 
Recall from \cite{FMMN} that $\Hq_s$ is given by 
\be
\nonumber
%\label{bks9}
\Hq_s:=
\left\{ \sigma_s = F e^{-\frac{s|Y|^2}{2\hbar_0}}\sqrt{\Omega_s}, \,\, F \,\,{\rm is} \,\,
J_{s} \text{-holomorphic and}\,\, ||\sigma_s||^{\Qu}_s<\infty \right\}, 
\ee
where the norm refers to the Hermitian structure on $\Hq$ as a
sub-bundle of $\Hpr$. 
The pairing (\ref{bks3}) between fibers of $\Hpr$ restricted to polarized sections 
defines the BKS pairing between fibers of $\Hq$, for which we use the same notation.

\begin{definition}
\label{defbks}
For $s,s'>0$, write 
\be
\nonumber
%\label{qsections}
\sigma_s = F e^{-\frac{s|Y|^2}{2\hbar_0}}\sqrt{\Omega_s}\in \Hq_s \text{ and } 
\sigma'_{s'} = F' e^{-\frac{s'|Y|^2}{2\hbar_0}}\sqrt{\Omega_{s'}}\in \Hq_{s'},
\ee
Their (quantum) BKS pairing is given by 
\be
\label{qbks}
\langle \sigma_s,\sigma'_{s'}\rangle^{\bks} = \int_{T^*K} \bar F F' 
 e^{-\frac{(s+s')}{2\hbar_0}|Y|^2}
\sqrt{\frac{\bar \Omega_s \wedge \Omega_{s'}}{b\epsilon}} \epsilon.
\ee
\end{definition}

As in (\ref{bks7}), the BKS pairing induces a connection, $\db$, on $\Hq$ defined
by the same formula
%\label{bks99}
\be
\nonumber
%\label{bks14}
\langle \sigma_s, \db_{\frac{\partial}{\partial s}}\sigma'_s \rangle^{\Qu} = 
\left.\frac{\partial}{\partial s'}\right|_{s'=s} \langle \sigma_s, \sigma'_{s'}\rangle^{\bks}.
\ee
for $\sigma_s \in \Hq_s$ and $\sigma'_{s'}\in \Hq_{s'}$.
Recall that, along the lines of \cite{AdPW},  
the quantum bundle is equipped with a quantum connection $\dq$, 
obtained from the orthogonal projection $P: \Hpr\rightarrow\Hq$. More precisely, we have 
$\dq = P \circ \dpr$, which implies immediately

\begin{theorem}
\label{corol}
The connections $\db$ and $\dq$ on $\Hq$ coincide. 
\end{theorem}

\begin{proof}
This is an obvious corollary of theorem \ref{bks8},
since the pairing on $\Hq$ is 
obtained by restriction from $\Hpr$.
\end{proof}

The BKS pairing formalism in this case provides results consistent with the approach of \cite{AdPW}. 
We will establish below the fact that the BKS pairing map for $\Hq$ is unitary and coincides
with the parallel transport of $\dq$.

\subsection{The quantum BKS pairing between K\"ahler polarizations
and the vertical polarization}
\la{ss31aq}

In \cite{FMMN}, we showed that the parallel transport associated to $\dq$
corresponds to Hall's CST, so that, in fact, $\dq$-parallel sections of 
$\Hq$ satisfy a heat equation on $K_\C$. This provides a better  understanding of 
some of the results of \cite{Ha3}. 

In \cite{Ha3}, Hall has also shown that, up to a constant, the CST  corresponds to the quantum 
BKS pairing map between the vertical real polarized and the K\"ahler polarized quantum Hilbert spaces. 
In this section, 
we establish in theorem \ref{bks12} an important identity for the
quantum BKS pairing between K\"ahler polarizations. This identity will help us understand better
the $s \rightarrow 0$ limit and, in section \ref{sss33}, the relation
between the quantum BKS pairing and the CST.
%%we further relate the BKS pairing on the quantum bundle $\Hq$ to the 
%%unitary parallel transport of the quantum connection $\dq$ and to the CST. 

Let now $\Delta$ be the invariant Laplacian on $K$ and
$\cal C$ denote the analytic continuation from $K$ to $K_\C$.
We recall from \cite{Ha1}, that the CST
is a unitary isomorphism of Hilbert spaces defined by 
\begin{eqnarray*}
C_\hbar : L^2(K,dx) & \rightarrow & {\cal H} L^2 (K_\C,d\nu_\hbar)\\
f & \mapsto & C_\hbar (f) 
= {\cal C} \circ e^{\frac{\hbar}{2}\Delta} f,
\end{eqnarray*}
where  $dx$ is the normalized Haar measure on $K$, 
$d\nu_{\hbar}(g)= \nu_{\hbar}(g)dg$ is the $K$-averaged heat kernel measure of \cite{Ha1},
and $dg$ is the Haar measure on $K_\C$. 
Recall from \cite{Ha2} the explicit form of the 
$K$-averaged heat kernel measure for $\hbar=s\hbar_0$
\be
\label{niu}
\nu_\hbar (g)= (a_s s^{n/2}\eta(Y))^{-1} e^{-\frac{|Y|^2}{\hbar}},
\ee
where
\be
\label{as}
a_s= (\pi \hbar_0)^{n/2}e^{|\rho|^2\hbar_0 s} ,
\ee
and
$\rho$ is the Weyl vector given by half the sum of the 
positive roots of $\K$.
Using the relation of parallel sections of
$\Hq$ with holomorphic functions on $K_\C$ we will now obtain an
explicit formula for the BKS pairing as an integral over $K_\C$. 
As in definition \ref{defbks} of the quantum BKS pairing, let 
\be
\label{parallel}
\sigma_s = F e^{-\frac{s|Y|^2}{2\hbar_0}}\sqrt{\Omega_s}\in \Hq_s \text{ and } 
\sigma'_{s'} = F' e^{-\frac{s'|Y|^2}{2\hbar_0}}\sqrt{\Omega_{s'}}\in \Hq_{s'}, s,s'>0,
\ee
where 
$F= \hat F \circ \psi_s$ and $F'= \hat F' \circ \psi_{s'}$ are, respectively, $J_s$-holomorphic
and $J_{s'}$-holomorphic functions on $T^*K$, with $\hat F, \hat F'$ 
homolorphic functions on $K_\C$. 

\begin{proposition}
\label{bks10}
The BKS pairing on $\Hq$ is given by
\be
\nonumber
%\label{bks11}
\langle \sigma_s,\sigma'_{s'}\rangle^{\bks} = a_{\frac{s+s'}{2}} \int_{K_\C} 
\overline{\hat F(gZ)} \hat F' (gZ^{-1})\,d\nu_{\hbar^{''}}(g),
\ee
where $g=xe^{iY}$, $Z = e^{i\frac{s-s'}{s+s'}Y}$ and $\hbar^{''} = \frac{s+s'}{2}\hbar_0$ .
\end{proposition}

\begin{proof}
{}From definition \ref{preqbks} and proposition \ref{bks4} we have
\be
\label{bks333}
\langle \sigma_s,\sigma'_{s'}\rangle^{\bks} = \int_{T^*K} \overline{\hat F(xe^{isY})} \hat F'(xe^{is'Y}) 
e^{-\frac{(s+s')}{2\hbar_0}|Y|^2} \left( \frac{s+s'}{2}\right)^{\frac{n}{2}} 
\eta\left(\frac{s+s'}{2}Y\right) \epsilon.
\ee
We recall from \cite{Ha2} and \cite{FMMN} that 
$\psi_s^* dg = |\Omega_s|^2\epsilon = s^{n}\eta^2(sY)\epsilon$.
Using (\ref{niu}) and (\ref{as}),
the formula then follows from a change of variables of integration, 
from $(x,Y)$ to $g=xe^{iY'}$ where
$Y' = \frac{s+s'}{2}Y$, and where at the end we renamed $Y'$ as $Y$ again.
\end{proof}

% {}Formula (\ref{bks333}), for $\hat F$, $\hat F'$ and $s$ fixed,  
% implies that the 
% BKS pairing between polarized sections of the form (\ref{parallel}) has a well defined $s'\rightarrow 0$ limit.
% Moreover, in this limit, setting $s=1$ and taking care of the different normalizations,
% we obtain the expression of \cite{Ha3} (see Theorem 2.6) for the 
% pairing between a vertically polarized section and a K\"ahler polarized one.
% (In our definition of the pairing there is an extra factor of $2^{-n/2}$ as 
% compared to \cite{Ha3}. See also the end of section \ref{ze1000}.) 
% This can also be checked directly from the formula 
% for the BKS pairing in the theorem \ref{bks12} below.

%%% \subsection{Relation between the BKS pairing and the CST}

%%We now recall that in \cite{FMMN} a relation was given between parallel sections 
%%of the quantum bundle and the coherent state transform for $K$. 
%%This will enable
%%us to find a relation between the BKS pairing and the CST.
% Let  now $\sigma_{s'}$ be in (\ref{parallel}), 
% have the form,
% $ \sigma_{s'}= \hat F' \circ \psi_{s'} e^{-\frac{s'|Y|^2}{2\hbar_0}}\sqrt{\Omega_{s'}}\in \Hq_{s'}$. 

The computation of the BKS pairing, and also its relation to the CST, is made very explicit 
with the following result.

\begin{theorem}
\label{bks12}
Let $\sigma_s \in \Hq_s$ and $\sigma'_{s'}\in \Hq_{s'}$ be given by
\ba
  \sigma_{s} &=& \left(C_{s\hbar_0} f \right) \circ \psi_{s} \ 
e^{-\frac{s|Y|^2}{2\hbar_0}}\sqrt{\Omega_{s}}   \nonumber \\
  \sigma'_{s'} &=& \left(C_{s'\hbar_0} f' \right) \circ \psi_{s'} \ 
e^{-\frac{s'|Y|^2}{2\hbar_0}}\sqrt{\Omega_{s'}} , \nonumber
\ea
where $s, s'>0$ and $f, f' \in L^2(K,dx)$.
We have,
\be
\nonumber
%\label{bks13}
\langle \sigma_s, \sigma'_{s'}\rangle^{\bks} = a_{\frac{s+s'}{2}} \langle f,f'\rangle_{L^2(K,dx)}. 
\ee
\end{theorem}

We remark that, from theorem \ref{corol}, one knows that the connection $\dq$ is also 
the connection induced from the BKS pairing
between infinitesimally close fibers of $\Hq$. Note, however, that of course one can give many 
different explicit 
Hermitian pairings between fibers of $\Hq$, such that they all induce the same connection $\dq$ 
on $\Hq$. 
Theorem \ref{bks12} implies that, among those, the pairing map that corresponds to the CST {\it is} the 
BKS pairing map naturally defined from geometric quantization. 
This is the map
$B^{\Qu}_{ss'}: \Hq_{s'}\rightarrow \Hq_s$ obtained by restriction of the
prequantum BKS map (\ref{BKSmap}),
\be
%\nonumber
\label{bksmapq}
\langle \sigma_s,\sigma'_{s'}\rangle^{\bks} = \langle \sigma_s, B^{\Qu}_{ss'}\sigma'_{s'}\rangle^{\Qu},
\ee
for $\sigma_s \in \Hq_s$ and $\sigma'_{s'} \in \Hq_{s'}$.

Postponing the proof of theorem \ref{bks12} to the end of this section, let us now address the problem of 
taking the $s'\rightarrow 0$ limit in the BKS pairing, which corresponds to degenerating one of 
the K\"ahler polarizations into the vertical polarization. 

{}From \cite{FMMN} we obtain, for $\sigma_{s}$ as in (\ref{parallel}),
$$
\sigma_{s} \in \Hq_{s} \ \Leftrightarrow \ \hat F \in {\cal H} L^2 (K_\C, d\nu_{\hbar})
$$
with $\hbar = s\hbar_0$.
Let us use the CST to define a family of continuous sections of $\Hq$ 
labelled by vectors $f \in L^2(K,dx)$,
\be
\la{ssp}
  \sigma_{s}= \left(C_{s\hbar_0} f \right) \circ \psi_{s} \ e^{-\frac{s|Y|^2}{2\hbar_0}}\sqrt{\Omega_{s}}.
\ee
Due to the fact that the coherent state transform is an isomorphism, 
the values of the sections $\sigma$ in (\ref{ssp}) span, for every $s>0$,
the whole fiber $\Hq_{s}$. This defines a global trivialization
of the vector bundle $\Hq$, which is particularly suitable
for studying the $s \rightarrow 0$ limit. 
Let $\Omega_0 = \lim_{s\rightarrow 0} \Omega_s = dx$, so that
$\sqrt{\Omega_0}= \sqrt{dx}$, and let $\Hq_0$ be
the space
$$
\Hq_0 = \left\{f \ \sqrt{\Omega_0}, \ f \in L^2(K,dx)\right\}.
$$
We consider the extension of the bundle
$\Hq \rightarrow \R_+$ to a bundle
\be
\nonumber
\widehat{\Hq} \ \rightarrow [0,\infty) ,  
\ee 
with a given trivialization
\be
\label{outono}
\begin{array}{rcl}
L^2(K,dx) \times [0,\infty) & \longrightarrow & \widehat \Hq 
 \\
(f, s) & \longmapsto &  \left(C_{s\hbar_0} f \right) \circ \psi_{s} \ e^{-\frac{s|Y|^2}{2\hbar_0}}\sqrt{\Omega_{s}} , \,s>0 \\
(f, 0) & \longmapsto & f\sqrt{\Omega_0}.
\end{array} 
%\nonumber
\ee

The sections (\ref{ssp}) are then extended  by continuity to $s=0$ by setting
\be
\la{ssp0} 
  \sigma_{0}= f\ \sqrt{\Omega_{0}} .
\ee

The BKS pairing between a K\"ahler polarized section $\sigma_s \in  \Hq_s$ as in (\ref{ssp}) and a 
vertically polarized section $\sigma'_{0}$ of the form (\ref{ssp0}),
with fixed 
$f'\in L^2 (K,dx)$, 
is then defined by
\be
\label{bksrk}
\langle \sigma_s,\sigma'_{0}\rangle^{\bks} = \int_{T^*K} \overline{\left(C_{s\hbar_0} f \right) (x e^{isY})}  
f'(x) \,
e^{-\frac{s}{2\hbar_0}|Y|^2} \left( \frac{s}{2}\right)^{\frac{n}{2}} 
\eta\left(\frac{s}{2}Y\right) \epsilon,
\ee
if the integral is absolutely convergent,
which corresponds formally to taking the $s'\rightarrow 0$ limit inside the integral in (\ref{bks333}). 

\begin{remark}
\la{rr1}
In \cite{Ha3}, it is shown that (after a trivial rescaling of the variable $Y$ to $sY$) 
this integral 
defines on dense subspaces a Hermitian pairing between the Hilbert space of 
K\"ahler polarized sections
and the Hilbert space of vertically polarized sections. 
Setting $s=1$ in (\ref{bksrk}) 
and taking care of the different norma\-li\-zations,
we obtain the same expression as in \cite{Ha3} (Theorem 2.6) for this pairing.
(In our definition of the pairing there is an extra factor of $2^{-n/2}$ as 
compared to \cite{Ha3}.) 
\end{remark}

Using our notation and normalization conventions, Hall's result implies that for any $s>0$,
\be
\label{cfff}
\langle\sigma_s,\sigma'_{0}\rangle^{\bks} = a_{\frac{s}{2}}\langle f,f'\rangle_{L^2(K,dx)}.
\ee

\begin{remark}
Notice that, even though the integral expression (\ref{bksrk}) can only be defined 
on dense subspaces, 
from (\ref{cfff}) it follows that
the pairing can be extended by continuity
to the respective completions \cite{Ha3}. This is in contrast to the pairing between two K\"ahler 
polarized spaces, where the integral defining the BKS pairing (\ref{bks3}) is absolutely convergent.
\end{remark}

We recall from  \cite{Ha3} that the pairing (\ref{bksrk}) 
does coincide with the general prescription for the geometric quantization
BKS pairing between a K\"ahler polarized Hilbert space and a vertically
polarized one. (See section 10.4 of \cite{Wo}.) This can also be easily deduced from (\ref{qbks}) and
proposition \ref{bks4}.

\begin{corollary}
\label{bks00}
Upon identifying $\Hq_{s'}$ with $L^2(K,dx)$ as in (\ref{outono}), 
the BKS pairing (\ref{qbks}) between polarized sections has a well defined (weak) $s'\rightarrow 0$ limit.
This limit coincides with (\ref{bksrk}).
\end{corollary}

\begin{proof}
This follows immediately from theorem \ref{bks12} by taking the limit  $s'\rightarrow 0$, 
with fixed $f,f'\in L^2(K,dx)$, which gives 
$\lim_{s'\rightarrow 0} \langle\sigma_s,\sigma'_{s'}\rangle^{\bks} = a_{\frac{s}{2}} \langle f, f'\rangle_{L^2(K,dx)}$.
as in (\ref{cfff})
\end{proof}

The associated pairing map between the Hilbert space of vertically polarized sections and the Hilbert space of
K\"ahler polarized sections (say for $s=1$) is unitary only up to scale as shown by Hall. 

%This can also be checked directly from the formula 
%for the BKS pairing in the next theorem.

To further clarify this issue, let us
address the choice of inner product for vertically polarized sections which
will be important for the unitarity of the BKS pairing maps $B^{\Qu}_{s0}$.
The inner product, $\langle\sigma_0, \sigma'_0\rangle^{\Qu}$,
 that is usually considered for vertically 
polarized sections of the form (\ref{ssp0}) reads 
$ \int_K \overline{f}f' \  dx$ (see e.g. \cite{Wo}). However, in order to 
have continuity of the Hermitian structure on 
$\widehat \Hq$ at $s=0$, we will define instead 
the inner product as
\be
\la{dis}
\langle\sigma_0, \sigma'_0\rangle^{\Qu} :=  \int_K \overline{f}f' \ (\pi \hbar_0)^{n/2}  dx .
\ee
We then have,
\begin{corollary}
\la{iii}
With the definition of the inner product (\ref{dis}),  the Hermitian structure
on the bundle $\widehat \Hq$ is continuous at  $s=0$.
\end{corollary}

\begin{proof}
{}The result follows from theorem \ref{bks12} with $s=s'$ and from (\ref{as}). 
\end{proof}
\noindent The reason for the appearance of the coefficient $(\pi \hbar_0)^{n/2}$ in (\ref{dis}) 
is related to the fact that 
in the inner product of the  holomorphic sections 
there is a factor associated to the contribution of the integration
along $\K$, which does not tend to $1$ when $s\rightarrow 0$. 
Corollary \ref{iii} is an indication 
that, for real polarizations which correspond to boundary points of spaces of holomorphic 
polarizations, and perhaps more generally, the inner product should be rescaled as in (\ref{dis}). 
As we will see, this choice, motivated by the continuity of the Hermitian structure, 
will also lead to the unitarity of the BKS pairing map, $B^{\Qu}_{s0}$ (rather than just
unitarity up to scale).

For the proof of theorem \ref{bks12}, we will need two auxiliary lemmas which prove 
analogues of equation (4) in \cite{Ha2}. Let 
$\rho_\hbar$ denote the analytic continuation of the heat kernel on $K$ to $K_\C$, as in \cite{Ha1}.
Moreover, let $*:K_\C\rightarrow K_\C$ denote the unique anti-holomorphic anti-automorphism of $K_\C$ 
which extends the map $x \rightarrow x^{-1}$ on $K$ (see also \cite{Ha4}).
We now show that the Dirac delta distribution on $K$, with respect to the Haar measure $dx$, 
can be informally written as
\be
\nonumber
%\label{bks16}
\delta(x) = \int_{K_\C}  \rho_{2\hbar} (x^{-1}g^* g)\,d\nu_\hbar(g).
\ee
More precisely, let $C(K)$ be the space of continuous functions on $K$ with the supremum
norm and let $\cal F$ denote the dense subspace of finite linear combinations 
of matrix elements of irreducible representations of $K$.
Then, we have
\begin{lemma}
\label{bks15}
Let $\{f_n \}_{n\in {\mathbb N}}$ be a sequence in $\cal F$ converging to $f \in C(K)$. 
Then
\be
\label{bksdelta}
\lim_{n\rightarrow \infty} 
\int_{K_\C} \left(\int_K \rho_{2\hbar} (x^{-1}g^* g)  f_n(x)) \ dx \right) d\nu_\hbar (g)
 = f(e),
\ee
where $e\in K$ is the identity.
\end{lemma}

\begin{proof}
Since evaluation at the identity is a 
continuous linear functional  on $C(K)$,
it suffices to show that for any matrix element $R_{ij}$ of any
irreducible representation $R$ of $K$, we have
\be
\label{carlos}
\int_{K_\C} \left(\int_K \rho_{2\hbar} (x^{-1}g^* g) R_{ij} (x) \ dx \right) d\nu_\hbar (g)
 = R_{ij}(e) = \delta_{ij}.
\ee

To show this, recall that $\rho_\hbar = \sum_R d_R e^{-\frac{\hbar}{2}c_R} \chi_R$, where 
$d_R$ is the dimension of $R$, 
$\chi_R$ denotes its character and 
$c_R$ is the negative of the eigenvalue of $\Delta$ 
corresponding to the eigenvector $\chi_R$.

Now, recall Weyl's classical orthogonality relations 
\be
\label{bks17}
\int_K \overline{R_{ij}(x)}R'_{lk} (x) dx= 
\frac{\delta_{RR'}}{d_R}\delta_{jk}\delta_{il}. 
\ee
{}From the unitarity of the 
CST of \cite{Ha1}, we then obtain the following identities,
\be
\label{bks17a}
e^{-\frac{\hbar}{2}(c_R + c_{R'})}
\int_{K_\C} 
\overline{R_{ij}(g)}R'_{lk}(g)      d\nu_\hbar (g)
= \frac{\delta_{RR'}}{d_R}\delta_{jk}\delta_{il}, 
\ee
where we have used the fact that $C_\hbar R_{ij} = e^{-\frac{\hbar}{2}c_R} R_{ij}$. 
(We are denoting by the same symbol a function on 
$K$ and its analytic continuation to $K_\C$.)

Upon substituting the expression for $\rho_{2\hbar}$ in (\ref{carlos}), 
using $\chi_R (g) = \sum_{i=1}^{d_R} R_{ii}(g)$ 
and the orthogonality relations (\ref{bks17}), 
(\ref{carlos}) becomes
$$
e^{-\hbar c_R}\sum_{k=1}^{d_R}\int_{K_\C} 
\overline{R_{ki}(g)} R_{kj}(g)d\nu_\hbar(g) = R_{ij}(e).  
$$
Using (\ref{bks17a}), we obtain the result.
\end{proof}

The Dirac delta distribution can also be written informally in a different way as
\be
\nonumber
\delta (x_1^{-1}x_2)=\int_{K_\C} 
\overline{\rho_\hbar(ge^{it Y}x_1^{-1})}
\rho_{\hbar'}(ge^{-it Y}x_2^{-1})\,d\nu_{\hbar''}(g),
\ee
which should be interpreted as above. More precisely,

\begin{lemma}
\label{onemore}
Let $\{f_n\}_{n\in {\mathbb N}}$ be a sequence in $\cal F$ converging to $f\in C(K)$. 
Then, 
\be
\label{bks19}
\lim_{n\rightarrow \infty}
\int_{K_\C} \left(\int_K \overline{\rho_\hbar(ge^{it Y}x_1^{-1})}
\rho_{\hbar'}(ge^{-it Y}x_2^{-1})  f_n(x_1) \ dx_1 \right) d\nu_{\hbar''} (g)
 = f(x_2),
\ee
for any $x_2 \in K$, $\hbar +\hbar' = 2\hbar''$, and any real number $t$.

\end{lemma}

\begin{proof}
To prove the lemma, we substitute the explicit expressions for the heat kernels in (\ref{bks19}) 
and use the fact that, with $g=xe^{iY}\in K_\C$, we have 
$\overline{\chi_R (xe^{i(1+t)Y}x_1^{-1})} = \chi_R (x_1e^{i(1+t)Y}x^{-1})$.
Using again $\rho_\hbar = \sum_R d_R e^{-\frac{\hbar}{2}c_R} \chi_R$ and 
rewriting the characters as a sum of products of matrix elements, we express the integral as an 
integral on $T^*K$. One then finds, upon integration along $K$, and use of the orthogonality
relations (\ref{bks17}), that the $t$ dependence 
in (\ref{bks19}) is apparent and cancels out. Finally, from the fact that 
$\hbar + \hbar' = 2\hbar''$ and using $\int_K dx =1$, one rewrites the expression as an 
integral on $K_\C$ so that equation (\ref{bks19}) becomes
\be
\nonumber
\lim_{n\rightarrow \infty}
\int_{K_\C} \left(\int_K \rho_{2\hbar} (x_2^{-1}x_1g^* g)  
f_n(x_1) \ dx_1 \right) d\nu_\hbar (g)
 = f(x_2),
\ee
which, by (\ref{bksdelta}), proves the lemma.
\end{proof}

We now prove theorem \ref{bks12}.

\begin{proof}({\it of theorem \ref{bks12}})
Let $\sigma_s,\sigma'_{s'}$ be as described in the theorem. From proposition \ref{bks10}, we have
\be
\label{bks18}
\langle \sigma_s,\sigma'_{s'}\rangle^{\bks} = a_{\frac{s+s'}{2}} \int_{K_\C} 
\overline{\hat F(gZ)} \hat F' (gZ^{-1})\,d\nu_{\hbar^{''}}(g) ,
\ee
where $Z = e^{i\frac{s-s'}{s+s'}Y}$, $\hbar^{''} = \frac{s+s'}{2}\hbar_0$ and where $\hat F = C_\hbar f$ and 
$\hat F'= C_{\hbar'} f'$, with $\hbar= s\hbar_0$ and
$\hbar' = s'\hbar_0$.  By definition, 
\be
\label{f1}
\hat F(gZ)=\int_K \rho_\hbar (gZx_1^{-1}) f(x_1) \ dx_1
\ee
and
\be
\label{f2}
\hat F'(gZ^{-1})=\int_K \rho_{\hbar'}(gZ^{-1}x_2^{-1})f'(x_2) \ dx_2.
\ee 
Substituting (\ref{f1}) and (\ref{f2}) in (\ref{bks18}), we see that to prove the result
for $f,f'\in {\cal F}$, it is enough to show that 
\be
\nonumber
\int_{K_\C} \left(\int_K \overline{\rho_\hbar(gZx_1^{-1})}
\rho_{\hbar'}(gZ^{-1}x_2^{-1}) f(x_1) \ dx_1 \right) d\nu_\hbar (g)
 = f(x_2),
\ee
which follows easily from lemma \ref{onemore}.
To prove the result for all $f,f'\in L^2(K,dx)$, we note that from (\ref{bks-1}) it follows that 
the prequantum BKS pairing is continuous on $\Hpr_s \times \Hpr_{s'}$ and, therefore, the quantum BKS pairing
on $\Hq_s \times \Hq_{s'}$ is also continuous. 
The theorem then follows from the fact that $\cal F$
is dense in $L^2(K,dx)$, which implies, from the isomorphisms 
given by (\ref{outono}), that the space of sections $\sigma_s$ of the form above with 
$f\in \cal F$, is also dense in $\Hq_s$ for all $s>0$.
\end{proof}

\subsection{Unitarity of the quantum BKS pairing map}
\label{ze1000}

In this section, we prove the unitarity of the quantum BKS pairing map (\ref{bksmapq}) 
between quantum Hilbert spaces $\Hq_s$ and $\Hq_{s'}$
using a direct calculation (the unitarity of these maps also follows
from theorem \ref{bks12}).

\begin{theorem}
\label{ze99}
The map $B^{\Qu}_{ss'}:\Hq_{s'}\rightarrow \Hq_s$ is unitary, for $s,s' \geq 0$.
\end{theorem}

In order to prove theorem \ref{ze99} we establish an auxiliary lemma.
Let
$$
\sigma^{R}_{s \,{ij}} = R_{ij}(xe^{isY})e^{-\frac{s|Y|^2}{2\hbar_0}} \sqrt{\Omega_s},
$$
where $R_{ij}$ denotes a matrix element of an irreducible representation $R$ of $K$, of dimension $d_R$, 
and $i,j =1,\dots,d_R$.

\begin{lemma}
\label{ze100}
{}For any $s,s' \geq 0$, we have 
\be
\nonumber
B^{\Qu}_{ss'}(\sigma^{R}_{s'\,{ij}}) = e^{-\frac{s-s'}{2}\hbar_0(c_R+|\rho|^2)} 
\sigma^{R}_{s\,{ij}}. 
\ee
\end{lemma}

\begin{proof}
The lemma is an immediate consequence of the equality 
\ba
\label{ze101}
\langle\sigma^{R}_{s\,{ij}}, B^{\Qu}_{ss'}\sigma^{R'}_{s'\,{i'j'}}\rangle^{\Qu} &=& 
\langle\sigma^{R}_{s \,{ij}}, \sigma^{R'}_{s' \,{i'j'}}\rangle^{\bks} = \\ \nonumber
&=& \frac{\delta_{RR'} 
\delta_{ii'} \delta_{jj'}}{d_R} (\pi \hbar_0)^{\frac{n}{2}} 
e^{\frac{s+s'}{2}\hbar_0(c_R+|\rho|^2)}.
\ea
%%This follows from theorem \ref{bks12} or from the following 
%%direct calculation, using $d\epsilon = dx\, dY$,
To establish it, recall that $\epsilon = dx\, dY$. For $s=s'=0$ the equality is trivial. 
Let either $s$ or $s'$ be positive. We then obtain
\ba
\nonumber
\langle\sigma^{R}_{s\,{ij}}, \sigma^{R'}_{s'\,{i'j'}}\rangle^{\bks} =
\int_{K\times \K} \sum_{k=1}^{d_R}
\overline{R_{ik}(x)R_{kj}(e^{isY})} \sum_{k'=1}^{d_{R'}}R'_{i'k'}(x)R'_{k'j'}(e^{is'Y})\\
\nonumber 
\textstyle e^{-\frac{(s+s')}{2\hbar_0}|Y|^2} 
\left(\frac{s+s'}{2}\right)^{\frac{n}{2}} \eta\left(\frac{s+s'}{2}Y
\right) dx \,dY =\\ \nonumber
=  \delta_{RR'} \delta_{ii'} 
\frac{1}{d_R} \int_\K \textstyle R_{jj'}(e^{i(s+s')Y})e^{-\frac{s+s'}{2\hbar_0}|Y|^2}
\left(\frac{s+s'}{2}\right)^{\frac{n}{2}} \eta\left(\frac{s+s'}{2}Y
\right) dY.
\ea
{}From (\ref{bks17a}) we obtain (\ref{ze101}) and the lemma.
\end{proof}

\begin{proof}({\it of theorem \ref{ze99}})
The sections $\{\sigma^{R}_{s' \,{ij}}\}$, with $R$ running over all 
irreducible representations of $K$ and 
$i,j=1,\dots,d_R$, form an orthogonal basis of $\Hq_{s'}$. Therefore, using  lemma \ref{ze100},
the unitarity of the map $B^{\Qu}_{ss'}$ is equivalent to the condition
\be
\label{ze399}
\nonumber
\frac{\langle B^{\Qu}_{ss'}\sigma^{R}_{s' \,{ij}},B^{\Qu}_{ss'}\sigma^{R}_{s' \,{ij}} \rangle^{\Qu}}
{\langle\sigma^{R}_{s' \,{ij}},\sigma^{R}_{s' \,{ij}} \rangle^{\Qu}} =
\frac{(\langle \sigma^{R}_{s \,{ij}}, B^{\Qu}_{ss'} \sigma^{R}_{s' \,{ij}}\rangle^{\Qu})^2}
{\langle \sigma^{R}_{s \,{ij}}, \sigma^{R}_{s \,{ij}}\rangle^{\Qu}
\langle \sigma^{R}_{s' \,{ij}}, \sigma^{R}_{s' \,{ij}}\rangle^{\Qu}} = 1.
\ee 

This follows directly from (\ref{ze101}).
\end{proof}

Note that the unitarity of $B^{\Qu}_{s0}$ follows immediately from
\cite{Ha3} if one chooses the inner product for the real polarization
as in (\ref{dis}).

\subsection{Relation between the quantum BKS pairing map and the CST}
\la{sss33}

In \cite{FMMN}, it was shown that the parallel transport between two fibers 
of $\Hq$, $\Hq_s$ and $\Hq_{s'}$, 
for the quantum connection $\dq$ corresponds to $C_\hbar \circ C_{\hbar'}^{-1}$,
where $\hbar = s \hbar_0$ and $\hbar' = s' \hbar_0$. 
More explicitly, 
we defined the CST bundle $\Hh \rightarrow \R_+$ with 
$\Hh_s = {\cal H}L^2(K_\C, d\nu_\hbar )$. 
The CST bundle comes equipped with a connection $\dh$ for which the parallel transport 
$U_{\hbar\hbar'}: \Hh_{s'} \rightarrow \Hh_{s}$ is given by the CST as 
$U_{\hbar\hbar'} = C_\hbar \circ C_{\hbar'}^{-1}$. 
Moreover, there exists a natural unitary bundle isomorphism
$S : \Hh \rightarrow \Hq$ such that $\dq = S \circ \dh \circ S^{-1}$. 
At the fiber over $s>0$ this is given by 
\be
\label{isom}
S_s (\hat F) = \hat F \circ \psi_s \,\frac{e^{-\frac{s|Y|^2}{2\hbar_0}}}{\sqrt{a_s}}\sqrt{\Omega_s},
\ee
for $\hat F \in {\cal H}L^2(K_\C, d\nu_\hbar)$.

As a corollary of the theorems above, 
we obtain the following result, which adds to theorem 4 of \cite{FMMN} the explanation of the role 
of the BKS pairing, from the point of view of the CST.

\begin{theorem}
\label{bks145}
{}For $s,s'>0$, the BKS pairing map $B^{\Qu}_{ss'}: \Hq_{s'}\rightarrow \Hq_s$ is 
given by 
\be
\nonumber
%\label{bks146}
B^{\Qu}_{ss'} = S_s \circ C_\hbar \circ C_{\hbar'}^{-1} \circ S_{s'}^{-1} 
= S_s \circ U_{\hbar\hbar'} \circ S_{s'}^{-1},
\ee
where $\hbar = s \hbar_0$ and $\hbar' = s' \hbar_0$. 
This pairing coincides with the parallel transport of $\dq$. 
\end{theorem}

\begin{proof}
The result follows from direct computation and from theorem \ref{bks12} and formula (\ref{isom}), 
where we use the unitarity of $C_\hbar$ and also $a_{\frac{s+s'}{2}}= \sqrt{a_s a_{s'}}$. 
Since $\dq \circ S = S \circ \dh$, the theorem states that the parallel transport 
of the quantum connection is given by the BKS pairing map which is clearly unitary.  
\end{proof}

This is in agreement with the results of \cite{FMMN} where, however, the quantum connection was  
defined in a different way by orthogonal projection. 
Recall also that sections of $\Hq$ which are horizontal with respect to $\dq$ satisfy a 
heat equation on $K_\C$. The corresponding explicit form for the connection $\dq$  
can also be found directly from theorem \ref{bks145}, using lemma 1 from \cite{FMMN}.

As a final remark, note that in the case when the group $K$ is abelian, 
all the formulas in the 
paper remain valid upon setting $\rho =0$ and $\eta \equiv 1$.

\vskip 2cm

%%%%%%%%%%%%%%%%%%%%%%%%%%%%%%%%%%%%%%%%%%%%%%%%%%%%%%%%%%%%%%%%%%%%
%%%%%%%%%%%%%%%%%%%%%%%%%%%%%%%%%%%%%%%%%%%%%%%%%%%%%%%%%%%%%%%%%%%%

\noindent {\large{\bf Acknowledgements:}} 
We wish to thank the referee for very useful comments, for suggesting
a more detailed analysis of the $s'\to 0$ limit and for pointing out a mistake 
in the proof of theorem 4.
The authors were partially
supported by the Center for Mathematics and its Applications, 
IST (CF, PM and JM), 
the Center for Analysis, Geometry and Dynamical Systems, IST (JPN), and
by the Funda\c c\~ao para a Ci\^encia e a Tecnologia (FCT) through
the programs PRAXIS XXI, POCTI, FEDER, the project POCTI/33943/MAT/2000, 
and the project CERN/FIS/43717/2001.

\vspace{1cm}

\noindent \small{$\dagger$Department of Mathematics, Instituto 
Superior T\'ecnico, Av. Rovisco Pais,
1049-001 Lisboa, Portugal \\ Email: cfloren, jmourao, jpnunes@math.ist.utl.pt}

\vskip 0.2cm

\noindent \small{$\ddagger$Center for Mathematics and its Applications, Instituto Superior
T\'ecnico, Av. Rovisco Pais, 1049-001 Lisboa, Portugal \\ Email:
pmatias@math.ist.utl.pt}


\begin{thebibliography}{MMM}


\small

\baselineskip 4.5mm

\bibitem[AdPW]{AdPW} S.Axelrod, S.Della Pietra, E.Witten, ``Geometric
quantization of
Chern-Simons gauge theory'', J. Diff. Geom. {\bf 33} (1991), 787-902.


\bibitem[FMMN]{FMMN} C.Florentino, P.Matias, J.Mour\~ao, J.P.Nunes, ``Geometric quantization, 
complex structures and the coherent state transform'', J. Funct. Anal. {\bf 221} (2005) 303--322.


\bibitem[FY]{FY} K.Furutani, S.Yoshizawa, ``A K\"ahler structure on the punctured cotangent 
bundle of complex and  quaternion projective spaces and its application to a 
geometric quantization.  II'',  Japan. J. Math. (N.S.)  {\bf 21}  (1995) 355--392.

\bibitem[Ha1]{Ha1} B.C.Hall, ``The Segal-Bargmann coherent state
transform for compact Lie groups'', J.
Funct. Anal. {\bf 122} (1994), 103-151.



\bibitem[Ha2]{Ha2} B.C.Hall, ``Phase space bounds for quantum mechanics on a 
 compact Lie group'',
 Comm. Math. Phys. {\bf 184} (1997) 233-250;



\bibitem[Ha3]{Ha3} B.C.Hall, ``Geometric quantization and the generalized 
Segal-Bargmann transform for
Lie groups of compact type'', Comm. Math. Physics {\bf 226} (2002) 233-268.

\bibitem[Ha4]{Ha4} B.C.Hall, ``The inverse Segal-Bargmann transform for compact Lie groups'', 
J. Funct. Anal. {\bf 143} (1997) 98-116.

\bibitem[Hi]{Hi} N.Hitchin, ``Flat connections and geometric quantization'', 
Comm, Math. Phys. {\bf 131} (1990) 347-380.


\bibitem[KW]{KW} W.Kirwin, S. Wu, ``Geometric 
quantization, parallel transport and the Fourier transform'', 
math.SG/0409555.

\bibitem[Ra]{Ra} J.Rawnsley, ``Coherent states and K\"ahler manifolds'', 
Quart. J. Math. Oxford Ser. (2) {\bf 28} (1977) 403-415; ``A 
nonunitary pairing of polarizations for 
the Kepler problem'' Trans. Am. Math. Soc. {\bf 250} (1970) 167-180.

\bibitem[Th]{Th} T.Thiemann, ``Reality Conditions Inducing Transforms 
for Quantum Gauge Field Theory
and Quantum Gravity'', Class. Quantum Grav. {\bf 13} (1996)
1383-1403.



\bibitem[Wo]{Wo} N. Woodhouse, ``Geometric Quantization'', Clarendon Press, Oxford, 1992.

\end{thebibliography}
\end{document}